\documentclass[sn-chicago]{sn-jnl}
\usepackage{dsfont}
\usepackage{color,graphics,graphicx}
\global\long\def\E{\mathbb{E}}
\newcommand{\KLEINO}{{\scriptstyle{\mathcal{O}}}}
\renewcommand{\P}{\mathbb{P}}

\newcommand{\N}{\mathds{N}}
\renewcommand{\R}{\mathds{R}}
\renewcommand\im{\mathrm{i}}

\renewcommand{\var}{\mathbb{V}\hspace*{-0.05cm}\textnormal{a\hspace*{0.02cm}r}}

\DeclareMathOperator{\rv}{RV_{\hspace*{-.05cm}n}}
\newcommand{\cov}{\mathbb{C}\textnormal{o\hspace*{0.02cm}v}}
\global\long\def\argmin{\arg\min}
\global\long\def\d{\mathrm{d}}
\newcommand{\abs}[1]{\ensuremath{\left\vert#1\right\vert}}
\newcommand{\tod}{~\stackrel{d}{\longrightarrow}~}

\definecolor{dblue}{rgb}{0.21,0.21,0.55}

\jyear{2022}%

\theoremstyle{thmstyleone}%
\newtheorem{theorem}{Theorem}
\newtheorem{cor}[theorem]{Corollary}

\theoremstyle{thmstyletwo}%
\newtheorem{example}{Example}%

\theoremstyle{thmstylethree}%
\newtheorem{assump}{Assumption}

\newcommand\EatDot[1]{.}
\newcommand\KillDot[1]{}

\begin{document}

\title[Efficient parameter estimation for parabolic SPDEs]{Efficient parameter estimation for parabolic SPDEs based on a log-linear model for realized volatilities}

\author*[1]{\fnm{Markus} \sur{Bibinger}}\email{markus.bibinger@mathematik.uni-wuerzburg.de}

\author*[1]{\fnm{Patrick} \sur{Bossert}}\email{patrick.bossert@mathematik.uni-wuerzburg.de}


\affil[1]{\orgdiv{Faculty of Mathematics and Computer Science}, \orgname{Julius-Maximilians-Universit\"at W\"urzburg}, \orgaddress{\street{Emil-Fischer-Str.\ 30}, \city{Würzburg}, \postcode{97074}, \country{Germany}}}

\abstract{We construct estimators for the parameters of a parabolic SPDE with one spatial dimension based on discrete observations of a solution in time and space on a bounded domain. We establish central limit theorems for a high-frequency asymptotic regime. The asymptotic variances are shown to be substantially smaller compared to existing estimation methods. Moreover, asymptotic confidence intervals are directly feasible. Our approach builds upon realized volatilities and their asymptotic illustration as response of a log-linear model with spatial explanatory variable. This yields efficient estimators based on realized volatilities with optimal rates of convergence and minimal variances. We demonstrate efficiency gains compared to previous estimation methods numerically and in Monte Carlo simulations.}


\maketitle
\clearpage

\section{Introduction}\label{sec1}
Dynamic models based on stochastic partial differential equations (SPDEs) are recently of great interest, in particular their calibration based on statistics, see, for instance, \cite{hambly}, \cite{fuglstad}, \cite{randolf2021} and \cite{randolf2022}. \citet{trabs}, \citet{cialenco} and \citet{chong} have independently of one another studied the parameter estimation for parabolic SPDEs based on power variation statistics of time increments when a solution of the SPDE is observed discretely in time and space. \citet{trabs} pointed out the relation of their estimators to realized volatilities which are well-known as key statistics for financial high-frequency data in econometrics. We develop estimators based on these realized volatilities which significantly improve upon the M-estimation from \cite{trabs}. Our new estimators attain smaller asymptotic variances, they are explicit functions of realized volatilities and we can readily provide asymptotic confidence intervals. Since generalized estimation approaches for small noise asymptotics in \citet{kaino2021}, rate-optimal estimation for more general observation schemes in \citet{hildebrandt}, long-span asymptotics in \citet{kaino}, and with two spatial dimensions in \citet{kainopre} have been built upon the M-estimator from \cite{trabs}, we expect that our new method is of interest to further improve parameter estimation for SPDEs.
Our theoretical framework is the same as in \cite{trabs}. We consider for $(t,y)\in\R_+\times[0,1]$ a linear parabolic SPDE 
\begin{align}
\d X_t(y)=\Big(\theta_2\frac{\partial^{2}X_t(y)}{\partial y^{2}}+\theta_1\frac{\partial X_t(y)}{\partial y}+\theta_0 X_t(y)\Big)\,\d t+\sigma\,\d B_{t}(y)\,,\label{eq:spde}
\end{align}
with one space dimension. The bounded spatial domain is the unit interval $[0,1]$, which can be easily generalized to some arbitrary bounded interval. Although estimation methods in case of an unbounded spatial domain are expected to be similar, the theory is significantly different, see \citet{trabs2}. $(B_{t}(y))$ is a cylindrical Brownian motion in a Sobolev space on $[0,1]$. The initial value $X_0(y)=\xi(y)$ is assumed to be independent of $(B_{t}(y))$. We work with Dirichlet boundary conditions: $X_{t}(0)=X_{t}(1)=0$, for all $t\in\R_+$. A specific example is the SPDE
\begin{align}
\d X_t(y)=\Big(\frac{\partial X_t(y)}{\partial y}+\frac{\kappa}{2}\frac{\partial^{2}X_t(y)}{\partial y^{2}}\Big)\,\d t+\sigma\,\d B_{t}(y),\label{eq:spde2}
\end{align}
used for the term structure model of \citet{cont2004}.

Existence and uniqueness of a mild solution of the SPDE \eqref{eq:spde} written \(\d X_t(y)=A_{\theta} X_t(y)\,\d t+\sigma\,\d B_{t}(y)\), with differential operator $A_{\theta}$, which is given by
\begin{align}\label{mild} X_t=\exp(t\,A_{\theta})\,\xi+\sigma\int_0^t \exp{((t-s)\,A_{\theta})}\,\d B_s\,,\end{align} 
with a Bochner integral and where $\exp(t\,A_{\theta})$ is the strongly continuous heat semigroup, is a classical result, see Chapter 6.5 in \cite{daPratoZabczyk1992}. We focus on parameter estimation based on \emph{discrete observations} of this solution $(X_t(y))$ on the unit square $(t,y)\in[0,1]\times[0,1]$. The spatial observation points $y_j$, $j=1,\ldots,m$, have at least distance $\delta>0$ from the boundaries at which the solution is zero by the Dirichlet conditions. 
\begin{assump}\label{assump-obs}
We assume equidistant \emph{high-frequency} observations in time $t_i=i\Delta_n$, $i=0,\dots,n$, where $\Delta_n=1/n\to 0$, asymptotically. We consider the same asymptotic regime as in \cite{trabs}, where $m=m_n\to\infty$, such that $m_n=\mathcal{O}(n^{\rho})$, for some $\rho\in(0,1/2)$, and $m\cdot\min_{j=2,\dots,m}\abs{y_j-y_{j-1}}$ is uniformly in $n$ bounded from below and from above. 
\end{assump}
This asymptotic regime with more observations in time than in space is natural for most applications. \citet{hildebrandt} and \cite{trabs} showed that in this regime the realized volatilities 
\[\rv(y_j)=\sum_{i=1}^n (X_{i\Delta_n}(y_j)-X_{(i-1)\Delta_n}(y_j))^2,~j=1,\ldots,m,\]
are sufficient to estimate the parameters with optimal rate $(m_n n)^{-1/2}$, while \cite{hildebrandt} establish different optimal convergence rates when the condition $m_n/\sqrt{n}\to 0$ is violated and propose rate-optimal estimators for this setup based on double increments in space and time. Our proofs and results could be generalized to non-equidistant observations in time, when their distances decay at the same order, but the observation schemes would affect asymptotic variances and thus complicate the results. Instead, there is no difference between equidistant and non-equidistant observations in space, since the spatial covariances will not be used for estimation and are asymptotically negligible for our results under Assumption \ref{assump-obs}.

The natural parameters, depending on $\theta_1\in\R$, and $\theta_2>0$, $\sigma>0$ from \eqref{eq:spde}, which are identifiable under high-frequency asymptotics are 
\begin{align}\label{parameters}\sigma^2_0:=\sigma^2/\sqrt{\theta_2}\quad\text{and}\quad \kappa:=\theta_1/\theta_2\,,\end{align}
the \emph{normalized volatility parameter} $\sigma^2_0$, and the \emph{curvature parameter} $\kappa$. The parameter $\theta_0\in\R$ could be estimated consistently only from observations on $[0,T]$, as $T\to\infty$. This is addressed in \cite{kaino}.

While \cite{trabs} focused first on estimating the volatility when the parameters $\theta_1$ and $\theta_2$ are known, we consider the estimation of the curvature parameter $\kappa$ in Section \ref{sec2}. We present an estimator for known $\sigma_0^2$ and a robustification for the case of unknown $\sigma_0^2$. In Section \ref{sec3} we develop a novel estimator for both parameters, $(\sigma_0^2,\kappa)$, which improves the M-estimator from Section 4 of \cite{trabs} significantly. It is based on a log-linear model for $\rv(y)$ with explanatory spatial variable $y$. Section \ref{sec4} is on the implementation and numerical results. We draw a numerical comparison of asymptotic variances and show the new estimators' improvement over existing methods. We demonstrate significant efficiency gains for finite-sample applications in a Monte Carlo simulation study. All proofs are given in Section \ref{sec5}.

\section{Curvature estimation}\label{sec2}
Section 3 of \cite{trabs} addressed the estimation of $\sigma^2$ in \eqref{eq:spde} when $\theta_1$ and $\theta_2$ are known. Here, we focus on the estimation of $\kappa$ from \eqref{parameters}, first when $\sigma_0^2$ is known and then for unknown volatility. The volatility estimator by \cite{trabs}, based on observations in one spatial point $y_j$, used the realized volatility $\rv(y_j)$. The central limit theorem with $\sqrt{n}$ rate for this estimator from Theorem 3.3 in \cite{trabs} yields equivalently that
\begin{align}\label{cltrv}\sqrt{n}\bigg(\frac{\rv(y_j)}{\sqrt{n}}-\frac{\exp(-\kappa y_j)\sigma_0^2}{\sqrt{\pi}}\bigg)\tod\mathcal{N}\big(0,\Gamma\sigma_0^4\exp(-2\kappa y_j)\big)\,,\end{align}
with $\Gamma\approx 0.75$ a numerical constant analytically determined by a series of covariances. Since the marginal processes of $X_t(y)$ in time have regularity 1/4, the scaling factor $1/\sqrt{n}$ for $\rv(y_j)$ is natural. To estimate $\kappa$ consistently we need observations in at least two distinct spatial points. A key observation in \cite{trabs} was that under Assumption \ref{assump-obs}, realized volatilities in different spatial observation points de-correlate asymptotically. From \eqref{cltrv} we can hence write
\begin{align}\label{rvmodel}
\frac{\rv(y_j)}{\sqrt{n }}=\exp(-\kappa y_j)\frac{\sigma_0^2}{\sqrt{\pi}}+\exp(-\kappa y_j)\,\sigma_0^2\,\sqrt{\frac{\Gamma}{n}}\,Z_j+R_{n,j}
\end{align}
with $Z_j$ i.i.d.\ standard normal and remainders $R_{n,j}$, which turn out to be asymptotically negligible for the asymptotic distribution of the estimators. The equation
\begin{align}\label{llmpre}
\log\Big(\frac{\rv(y_j)}{\sqrt{n }}\Big)&=-\kappa y_j+\log\Big(\frac{\sigma_0^2}{\sqrt{\pi}}\Big)+\log\big(1+\sqrt{\Gamma\pi \Delta_n}\,Z_j\big)\\
&\notag\hspace*{2.2cm} +\log\Big(1+ \frac{ R_{n,j}\exp(\kappa y_j)\sqrt{\pi}\sigma_0^{-2}}{1+\sqrt{\Gamma\pi \Delta_n}\,Z_j}\Big)\,,
\end{align}
and an expansion of the logarithm yield an approximation
\begin{align}\label{kappaapprox}\kappa\approx \frac{-\log\big(\Delta_n^{1/2}\rv(y_j)\big)+\log(\sigma_0^2)-\log(\sqrt{\pi})}{y_j}+\frac{\sqrt{\Gamma\pi \Delta_n}}{y_j}\,Z_j\,.\end{align}
In the upcoming example we briefly discuss optimal estimation in a related simple statistical model.
\begin{example}\label{ex1}
Assume independent observations $Y_i\sim\mathcal{N}(\mu,\varsigma_i^2),i=1,\ldots,m$, where $\mu$ is unknown and $\varsigma_i^2>0$ are known. The maximum likelihood estimator (mle) is given by 
\begin{align}
\hat{\mu}=\frac{\sum_{i=1}^m Y_i\varsigma_i^{-2}}{\sum_{i=1}^m\varsigma_i^{-2}}\,.\label{equation_mle}\end{align}
The expected value and variance of this mle are
\begin{align*}
\E[\hat{\mu}]= \mu ~~~~~\text{and}~~~~~ \var(\hat{\mu})=\bigg(\sum_{i=1}^m\varsigma_i^{-2}\bigg)^{-1}.
\end{align*}
Note that $\varsigma_i^{-2}$ can be viewed as Fisher information of observing $Y_i$. The efficiency of the mle in this model is implied by standard asymptotic statistics.
\end{example}
If we have independent observations with the same expectation and variances as in Example \ref{ex1}, but not necessarily normally distributed, the estimator $\hat\mu$ from \eqref{equation_mle} can be shown to be the linear unbiased estimator with minimal variance.

If $\sigma_0^2$ is known, this and \eqref{kappaapprox} motivate the following curvature estimator:
\begin{align}\label{kappahat}
\hat{\kappa}_{n,m} =\frac{-\sum_{j=1}^m \log\Big(\frac{\rv(y_j)}{\sqrt{n}}\Big)y_j+\sum_{j=1}^m\log\Big(\frac{\sigma_0^2}{\sqrt{\pi}}\Big)y_j}{\sum_{j=1}^my_j^2}\,.
\end{align}
\begin{theorem}
Grant Assumptions \ref{assump-obs} and \ref{cond} with $y_1=\delta$, $y_m=1-\delta$, and $\delta\in(0,1/2)$. Then, the estimator \eqref{kappahat} satisfies, as $n\to\infty$, the central limit theorem (clt)
\begin{align}\label{clthatkappa}
\sqrt{nm_n}\big(\hat{\kappa}_{n,m}-\kappa\big)\tod\mathcal{N}\bigg(0,\frac{3\Gamma\pi}{1-\delta+\delta^2}\bigg)\,.
\end{align}
\end{theorem}
Typically $\delta$ will be small and the asymptotic variance close to $3\Gamma\pi$. Assumption \ref{cond} poses a mild restriction on the initial condition $\xi$ and is stated at the beginning of Section \ref{sec5}. The logarithm yields a variance stabilizing transformation for \eqref{cltrv} and the delta-method readily a clt for log realized volatilities with constant asymptotic variances. This implies a clt for the estimator as $\Delta_n\to 0$, and when $1<m<\infty$ is fix. The proof of \eqref{clthatkappa} is, however, not obvious and based on an application of a clt for weakly dependent triangular arrays by \cite{utev}.

If $\sigma_0^2$ is unknown the estimator \eqref{kappahat} is infeasible. Considering differences for different spatial points in \eqref{llmpre} yields a natural generalization of Example \ref{ex1} and the estimator \eqref{kappahat} for this case:
\begin{align}\label{kappahatf}
\hat{\varkappa}_{n,m}=\frac{\sum_{j\neq l}\log\Big(\frac{\rv(y_j)}{\rv(y_l)}\Big)(y_l-y_j)}{\sum_{j\neq l}(y_j-y_l)^2}\,.
\end{align}
This estimator achieves as well the parametric rate of convergence $\sqrt{nm_n}$, it is asymptotically unbiased and satisfies a clt. Its asymptotic variance is, however, much larger than the one in \eqref{clthatkappa}.
\begin{theorem}
Grant Assumptions \ref{assump-obs} and \ref{cond} with $y_1=\delta$, $y_m=1-\delta$, and $\delta\in(0,1/2)$. Then, the estimator \eqref{kappahatf} satisfies, as $n\to\infty$, the clt
\begin{align}\label{clthatkappaf}
\sqrt{nm_n}\big(\hat{\varkappa}_{n,m}-\kappa\big)\tod\mathcal{N}\bigg(0,\frac{12\Gamma\pi}{(1-2\delta)^2}\bigg)\,.
\end{align}
\end{theorem}
In particular, consistency of the estimator holds as $n\to\infty$, also if $m\ge 2$ remains fix. The clts \eqref{clthatkappa} and \eqref{clthatkappaf} are feasible, i.e.\ the asymptotic variances are known constants and do not hinge on any unknown parameters. Hence, asymptotic confidence intervals can be constructed based on the theorems.
\section{Asymptotic log-linear model for realized volatilities and least squares estimation}\label{sec3}
Applying the logarithm to \eqref{rvmodel} and a first-order Taylor expansion 
\[\log(a+x)=\log(a)+\frac{x}{a}+\mathcal{O}\Big(\frac{x^2}{a^2}\Big)\,,~x\to 0\,,\]
yield an asymptotic \emph{log-linear model}
\begin{align}\label{llm}\tag{LLM}
\log\Big(\frac{\rv(y_j)}{\sqrt{n }}\Big)=-\kappa y_j+\log\Big(\frac{\sigma_0^2}{\sqrt{\pi}}\Big)+\sqrt{\frac{\Gamma\pi}{n}}\,Z_j+\tilde R_{n,j}
\end{align}
for the rescaled realized volatilities, with $Z_j$ i.i.d.\ standard normal and remainders $\tilde R_{n,j}$, which turn out to be asymptotically negligible for the asymptotic distribution of the estimators. When we ignore the remainders $\tilde R_{n,j}$, the estimation of $-\kappa$ is then directly equivalent to estimating the slope parameter in a simple ordinary linear regression model with normal errors. The intercept parameter in the model \eqref{llm} is a strictly monotone transformation 
\begin{align}\label{transf}\alpha(\sigma_0^2)=\log\big(\frac{\sigma_0^2}{\sqrt{\pi}}\big)\end{align}
of $\sigma_0^2$. To exploit the analogy of \eqref{llm} to a log-linear model, it is useful to recall some standard results on least squares estimation for linear regression.
\begin{example}\label{regexample}
In a simple linear ordinary regression model
\[Y_i=\alpha+\beta x_i+\epsilon_i~,~i=1,\ldots,m\,,\]
with white noise $\epsilon_i$, homoscedastic with variance $\var(\epsilon_i)=\sigma^2$, the least squares estimation yields
\begin{subequations}
\begin{align}\label{reg1}
\hat\beta&=\frac{\sum_{j=1}^m(x_j-\bar x)(Y_j-\bar Y)}{\sum_{j=1}^m(x_j-\bar x)^2}\,,\\
\label{reg2}\hat\alpha&=\bar Y-\hat\beta\bar x\,,
\end{align}
with the sample averages $\bar Y=m^{-1}\sum_{j=1}^m Y_j$, and $\bar x=m^{-1}\sum_{j=1}^m x_j$. The estimators \eqref{reg1} and \eqref{reg2} are known to be BLUE (best linear unbiased estimators) by the famous Gau\ss-Markov theorem, i.e.\ they have minimal variances among all linear and unbiased estimators. In the normal linear model, if $\epsilon_i\stackrel{i.i.d.}{\sim}\mathcal{N}(0,\sigma^2)$, the least squares estimator coincides with the mle and standard results imply asymptotic efficiency. The variance-covariance matrix of $(\hat\alpha,\hat\beta)$ is well-known and
\begin{align}
\label{lm1}
\var(\hat\beta)&=\frac{\sigma^2}{\sum_{j=1}^m (x_j-\bar x)^2}\,,\\
\label{lm2}
\var(\hat\alpha)&=\frac{\sigma^2\sum_{j=1}^m x_j^2}{m\sum_{j=1}^m (x_j-\bar x)^2}\,,\\
\label{lm3}
\cov(\hat\alpha,\hat\beta)&=-\frac{\sigma^2\bar x}{\sum_{j=1}^m (x_j-\bar x)^2}\,.
\end{align} 
\end{subequations}
For the derivation of \eqref{reg1}-\eqref{lm3} in this example see, for instance, Example 7.2-1 of \cite{zimmerman2020}.
\end{example}
We give this elementary example, since our estimator and the asymptotic variance-covariance matrix of our estimator will be in line with the translation of the example to our model \eqref{llm}.

The M-estimation of \cite{trabs} was based on the parametric regression model
\begin{align}\label{prm}
\frac{\rv(y_j)}{\sqrt{n }}=\frac{\sigma_0^2}{\sqrt{\pi}}\exp(-\kappa y_j)+\delta_{n,j}\,,
\end{align}
with non-standard observation errors $(\delta_{n,j})$. The proposed estimator
\begin{equation}\label{bt}\argmin_{s,k}\sum_{j=1}^{m}\bigg(\frac{\rv(y_j)}{\sqrt{n }}-\frac{s^2}{\sqrt{\pi}}\exp(-k y_j)\bigg)^{2} 
\end{equation}
was shown to be rate-optimal and asymptotically normally distributed in Theorem 4.2 of \cite{trabs}. In view of the analogy of \eqref{llm} to an ordinary linear regression model, however, it appears clear that the estimation method by \cite{trabs} is inefficient, since \emph{ordinary} least squares is applied to a model with \emph{heteroscedastic} errors. In fact, \emph{generalized} least squares could render a more efficient estimator related to our new methods. In model \eqref{prm}, the variances of $\delta_{n,j}$ depend on $j$ via the factor $\exp(-2\kappa y_j)$. This induces, moreover, that the asymptotic variance-covariance matrix of the estimator \eqref{bt} depends on the parameter $(\sigma_0^2,\kappa)$. In line with the least squares estimator from Example \ref{regexample}, the asymptotic distribution of our estimator will \emph{not} depend on the parameter. 

Writing $\overline{y}=m_n^{-1}\sum_{j=1}^{m_n}y_j$, our estimator for $\kappa$ reads
\begin{subequations}
\begin{align}\label{lskappa}
\hat{\kappa}_{n,m}^{LS}&=\frac{\sum_{j=1}^{m_n}\log\big(\frac{\rv(y_j)}{\sqrt{n}}\big) y_j-\overline{y}\sum_{j=1}^{m_n}\log\big(\frac{\rv(y_j)}{\sqrt{n}}\big)}{m_n(\overline{y})^2-\sum_{j=1}^{m_n}y_j^2}\\
&=-\frac{\sum_{j=1}^{m_n}\bigg(\log\big(\frac{\rv(y_j)}{\sqrt{n}}\big)-\Big(m_n^{-1}\sum_{u=1}^{m_n}\log\big(\frac{\rv(y_u)}{\sqrt{n}}\big)\Big)\bigg) \big(y_j-\overline{y}\big)}{\sum_{j=1}^{m_n}\big(y_j-\overline{y}\big)^2}\,.\notag
\end{align}
The estimator for the intercept is
\begin{align}\label{lsint}
\widehat\alpha^{LS}(\sigma_{0}^{2})&=\overline{y}\hat{\kappa}_{n,m}^{LS}+m_n^{-1}\sum _{j=1}^{m_n} \log\Big(\frac{\rv(y_j)}{\sqrt{n}}\Big)\\
&=\frac{\overline{y} \Big(\sum _{j=1}^{m_n} \log\big(\frac{\rv(y_j)}{\sqrt{n}}\big) y_j\Big)-m_n^{-1}\Big(\sum _{j=1}^{m_n} \log\big(\frac{\rv(y_j)}{\sqrt{n}}\big)\Big)\big( \sum _{j=1}^{m_n} y_j^2\big)}{m_n(\overline{y})^2- \sum _{j=1}^{m_n} y_j^2}\,.\notag
\end{align}
\end{subequations}
We shall prove that the OLS-estimator \eqref{lskappa} for $\kappa$ in our log-linear model is in fact identical to the estimator $\hat{\varkappa}_{n,m}$ from \eqref{kappahatf}.
\begin{theorem}\label{jointclt}
Grant Assumptions \ref{assump-obs} and \ref{cond} with $y_1=\delta$, $y_m=1-\delta$, and $\delta\in(0,1/2)$. The estimators \eqref{lskappa} and \eqref{lsint} satisfy, as $n\to\infty$, the bivariate clt\\[-.2125cm]
\begin{align*}
\sqrt{nm_n}\left(
\begin{pmatrix}
\hat{\kappa}_{n,m}^{LS} \\[.2cm] \widehat\alpha^{LS}(\sigma_{0}^{2})
\end{pmatrix}
-
\begin{pmatrix}
\kappa\\[.2cm]  \alpha(\sigma_0^2)
\end{pmatrix}
\right)\tod\mathcal{N}\left(0,\Sigma\right),\end{align*}
with the asymptotic variance-covariance matrix\\[-.2125cm] 
\begin{align*}
\Sigma= \begin{pmatrix}\frac{12\Gamma\pi}{(1-2 \delta)^2 }&\frac{6 \Gamma\pi}{(1-2\delta)^2}\\ \frac{6 \Gamma\pi}{(1-2\delta)^2}&4\Gamma\pi \frac{1-\delta+\delta^2}{(1-2\delta)^2}\end{pmatrix}.
\end{align*}
\end{theorem}
In particular, consistency of the estimators holds as $n\to\infty$, also if $m\ge 2$ remains fix. Different to the typical situation with an unknown noise variance in Example \ref{regexample}, the noise variance in \eqref{llm} is a known constant. Therefore, different to Theorem 4.2 of \cite{trabs}, our central limit theorem is readily feasible and provides asymptotic confidence intervals.

An application of the multivariate delta method yields the bivariate clt for the estimation errors of the two point estimators.
\begin{cor}
Under the assumptions of Theorem \ref{jointclt}, it holds that\\[-.2125cm]
\begin{align*}
\sqrt{nm_n}\left(
\begin{pmatrix}
\hat{\kappa}_{n,m}^{LS} \\ (\hat\sigma_{0}^{2})^{LS} 
\end{pmatrix}
-
\begin{pmatrix}
 \kappa \\ \sigma_0^2
\end{pmatrix}\right)\tod\mathcal{N}\left(0,\tilde\Sigma\right)\,,
\end{align*}

where $(\hat\sigma_{0}^{2})^{LS}$ is obtained from $\widehat\alpha^{LS}(\sigma_{0}^{2})$ with the inverse of \eqref{transf}, with\\[-.2125cm]
\begin{align*}
\tilde\Sigma=\begin{pmatrix}
\frac{12\Gamma\pi}{(1-2\delta)^2}
&\frac{6\sigma_0^2\Gamma\pi}{(1-2\delta)^2}
\\ \frac{6\sigma_0^2\Gamma\pi}{(1-2\delta)^2}&
\frac{4\sigma_0^4\Gamma\pi(1-\delta+\delta^2)}{(1-2\delta)^2}
\end{pmatrix}.
\end{align*}
\end{cor}
Here, naturally the parameter $\sigma_0^2$ occurs in the asymptotic variance of the estimated volatility and in the asymptotic covariance. The transformation or plug-in, however, still readily yield asymptotic confidence intervals.


\section{Numerical illustration and simulations}\label{sec4}

\subsection{Numerical comparison of asymptotic variances}
\begin{figure}[t]
\centering
\includegraphics[width=11cm]{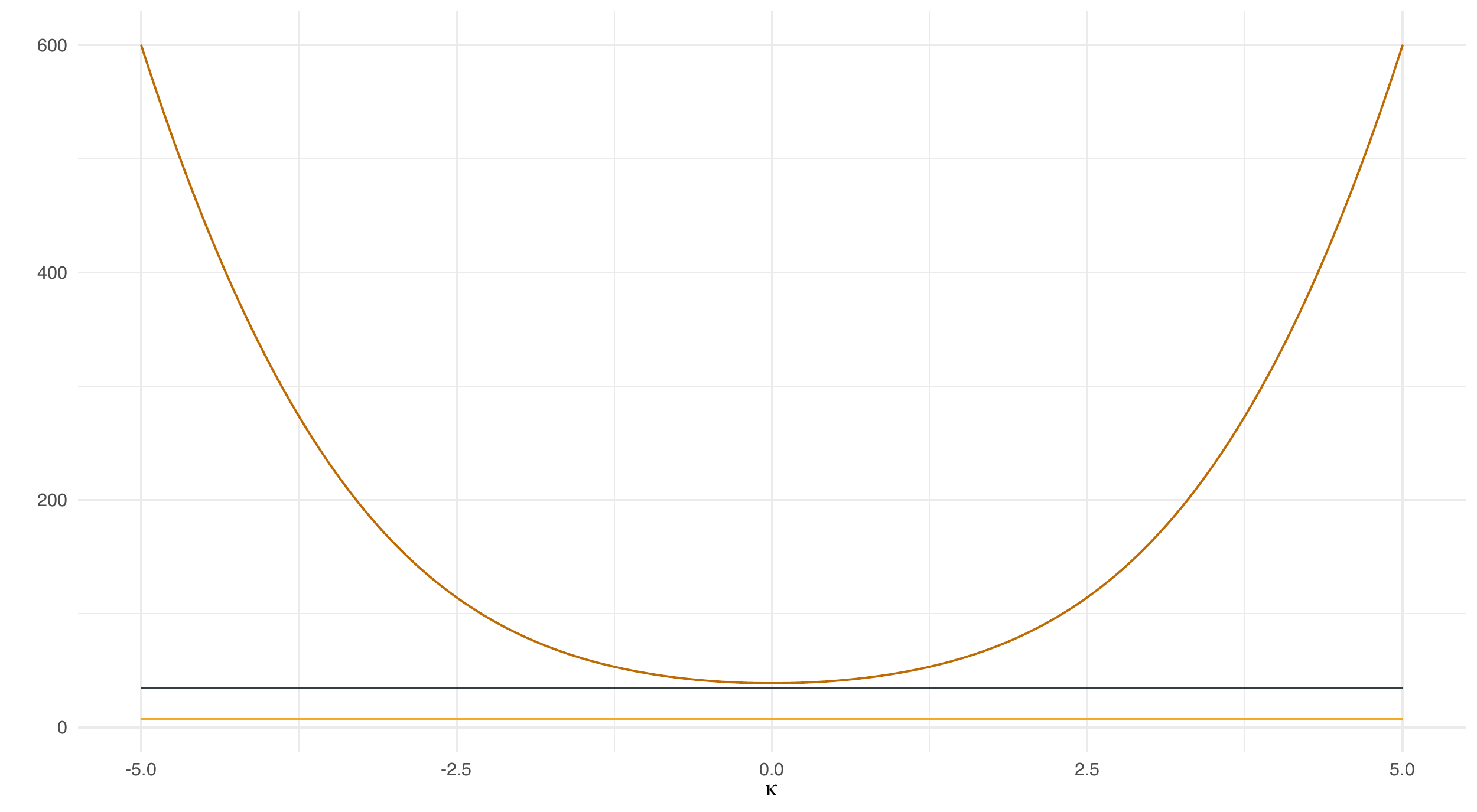}
\includegraphics[width=5.75cm]{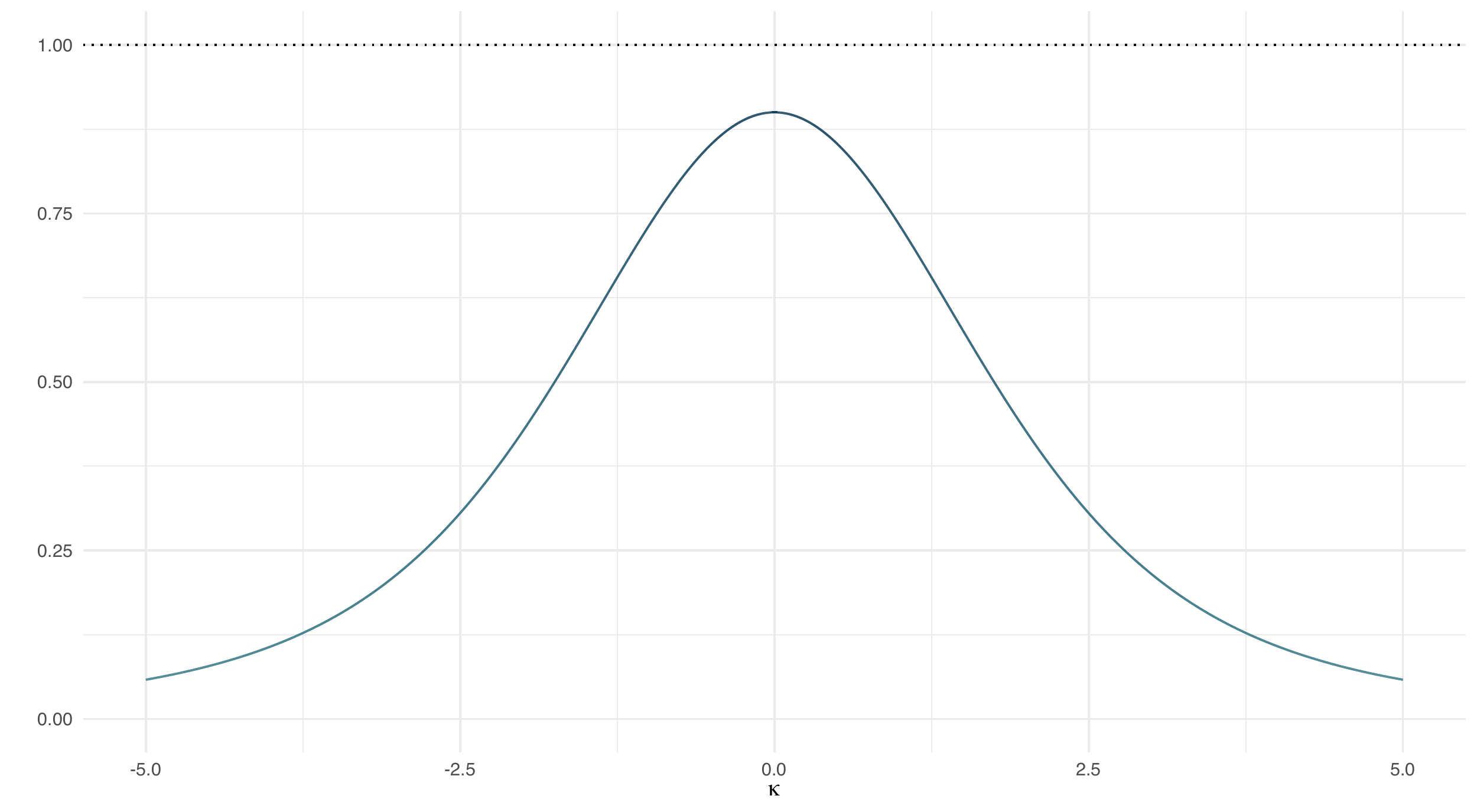}\includegraphics[width=5.75cm]{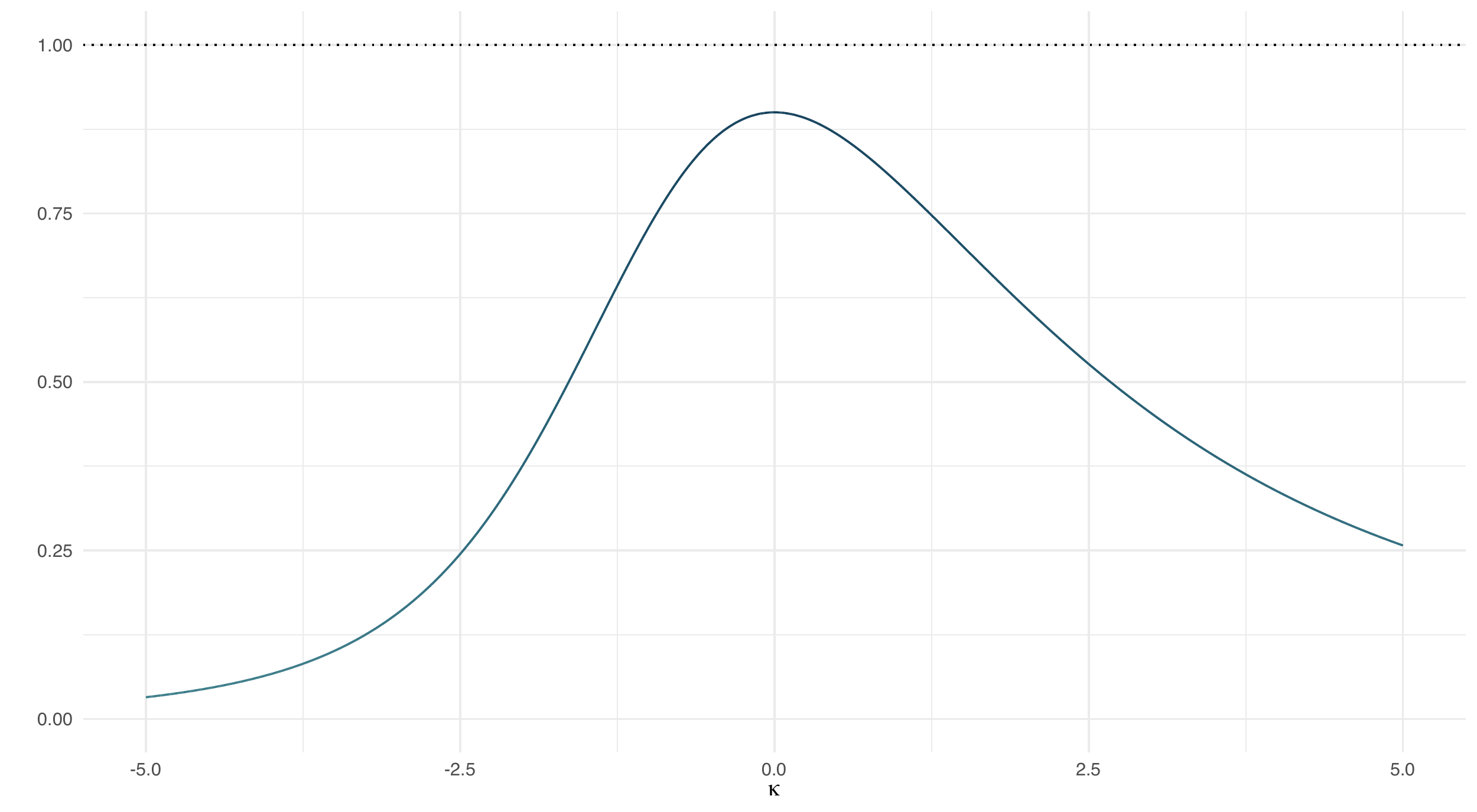}
\caption{Top panel: Comparison of asymptotic variances of $\hat{\kappa}_{n,m}$ from \eqref{kappahat} (for known $\sigma_0^2$), $\hat{\kappa}_{n,m}^{LS}$ from \eqref{lskappa} and the estimator from \cite{trabs}, for $\delta=0{.}05$, and for different values of $\kappa$. Lower panel: Ratio of asymptotic variances of new method using \eqref{lskappa} and \eqref{lsint} vs.\ \cite{trabs}, left for estimating $\kappa$, right for $\sigma_0^2$.}\label{fig1}
\end{figure}

The top panel of Figure \ref{fig1} gives a comparison of the asymptotic variances for curvature estimation, $\kappa$, of our new estimators to the minimum contrast estimator of \cite{trabs}. We fix $\delta=0{.}05$. While the asymptotic variance-covariance matrix of our new estimator is rather simple and explicit, the one in \cite{trabs} is more complicated but can be explicitly computed from their Eq.\ (21)-(23).

The uniformly smallest variance is the one of $\hat{\kappa}_{n,m}$ from \eqref{kappahat}. It is visualized with the yellow curve which is constant in $\kappa$, i.e.\ the asymptotic variance does not hinge on the parameter. This estimator, however, requires that the volatility $\sigma_0^2$ is known. It is thus fair to compare the asymptotic variance of the minimum contrast estimator from \cite{trabs} only to the least squares estimator based on the log-linear model, since both work for unknown $\sigma_0^2$. While the asymptotic variance of the new estimator, visualized with the black curve, does not hinge on the parameter value, the variance of the estimator by \cite{trabs} (brown curve) is in particular large when $\kappa$ has a larger distance to zero. All curves in the top panel of Figure \ref{fig1} do not depend on the value of $\sigma_0^2$. Our least squares estimator in the log-linear model uniformly dominates the estimator from \cite{trabs}. For $\delta\to 0$, the asymptotic variances of the two least squares estimators would coincide in $\kappa=0$. However, due to the different dependence on $\delta$, the asymptotic variance of the estimator from \cite{trabs} is larger than the one of our new estimator also in $\kappa=0$. The lower panel of Figure \ref{fig1} shows the ratios of the asymptotic variances of the two least squares estimators for both unknown parameters. Left we see the ratio for curvature estimation determined by the ratio of the black and brown curves from the graphic in the top panel. Right we see the ratio of the asymptotic variances for estimating $\sigma_0^2$, as a function depending on different values of $\kappa$. This ratio does not hinge on the value of $\sigma_0^2$.

\subsection{Monte Carlo simulation study}
The simulation of the SPDE is based on its spectral decomposition \eqref{eq:Representation} and an exact simulation of the Ornstein-Uhlenbeck coordinate processes. In \cite{trabs} a truncation method was suggested to approximate the infinite series $\sum_{k=1}^{\infty}x_k(t)e_k(y)$ in \eqref{eq:Representation} by a finite sum $\sum_{k=1}^{K}x_k(t)e_k(y)$, up to some spectral cut-off frequency $K$ which needs to be set sufficiently large. In \citet{kaino} this procedure was adopted, but they observed that choosing $K$ too small results in a strong systematic bias of simulated estimates. A sufficiently large $K$ depends on the number of observations, but even for moderate sample sizes $K=10^5$ was recommended by \cite{kaino}. This leads to tedious, long computation times as reported in \cite{kaino}. A nice improvement for observations on an equidistant grid in time and space has been presented by \citet{hildebrandt2020} using a replacement method instead of the truncation method. The replacement method approximates addends with large frequencies in the Fourier series using a suitable set of independent random vectors instead of simply cutting off these terms. The spectral frequency to start with replacement can be set much smaller than the cut-off $K$ for truncation. We thus use Algorithm 3.2 from \cite{hildebrandt2020} here, which allows to simulate a solution of the SPDE with the same precision as the truncation method while reducing the computation time considerably. For instance, for $n=10000$ and $m=100$, the computation time with a standard computer of the truncation method is almost 6 hours while the replacement method requires less than one minute. In \cite{hildebrandt2020} we also find bounds for the total variation distance between approximated and true distribution allowing to select an appropriate trade-off between precision and computation time. We implement the method with $20\cdot m_n$ as the spectral frequency to start with replacement. We simulate observations on equidistant grid points in time and space. We illustrate results for a spatial resolution with $m=11$, and a temporal resolution with $n=1000$. This is in line with Assumption \ref{assump-obs}. We simulated results for $m=100$ and $n=10000$, as well. Although the ratio of spatial and temporal observations is more critical then, the normalized estimation results were similar. If the condition $m_n\,n^{-1/2}\to 0$ is violated, however, we see that the variances of the estimators decrease at a slower rate than $(n\cdot m_n)^{-1/2}$. While the numerical computation of estimators as in \cite{trabs} relies on optimization algorithms, the implementation of our estimators is simple, since they rely on explicit transformations of the data. We use the programming language R and provide a package for these simulations on github.\,\footnote{\href{https://github.com/pabolang/ParabolicSPDEs}{https://github.com/pabolang/ParabolicSPDEs}}

\begin{figure}[t]
\centering
\includegraphics[width=11.9cm]{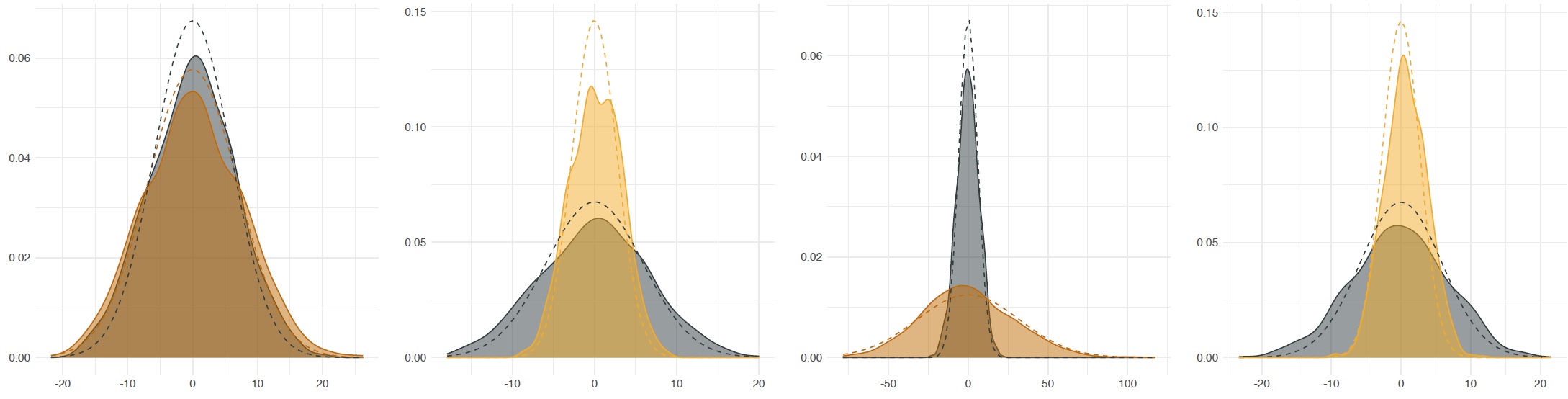}
\caption{Comparison of empirical distributions of normalized estimation errors for $\kappa$ from simulation with $n=1000$, $m=11$, $\sigma_0^2=1$, $\kappa=1$ in the left two columns and $\kappa=6$ in the right two columns. Grey is for $\hat{\kappa}_{n,m}^{LS}$, brown for the estimator by \cite{trabs} and yellow for  $\hat{\kappa}_{n,m}$.}\label{fig2}
\end{figure}

Figure \ref{fig2} compares empirical distributions of normalized estimation errors 
\begin{enumerate}
\item of $\hat{\kappa}_{n,m}^{LS}$ (grey) vs.\ \cite{trabs} (brown) and 
\item of $\hat{\kappa}_{n,m}$ with known $\sigma_0^2$ (yellow) compared to $\hat{\kappa}_{n,m}^{LS}$ (grey),
\end{enumerate}
for small curvature $\kappa=1$ in the left two columns, and larger curvature $\kappa=6$ in the right two columns. The plots are based on a Monte Carlo simulation with 1000 iterations, and for $n=1000$, $m=11$, and $\sigma_0^2=1$. We use the standard R density plots with Gaussian kernels and bandwidths selected by Silverman's rule of thumb. The dotted lines give the corresponding densities of the asymptotic limit distributions. We can report that analogous plots for different parameter values of $\sigma_0^2$ look (almost) identical. With increasing values of $n$ and $m$, the fit of the asymptotic distributions becomes more accurate, otherwise the plots look as well similar as long as $m\le \sqrt{n}$.

As expected, the efficiency gains of the new method are much more relevant for larger curvature. In particular, in the third plot from the left for $\kappa=6$, the new estimator outperforms the one from \cite{trabs} significantly. In the first plot from the left for $\kappa=1$, instead, the two estimators have similar empirical distributions. The fit of the asymptotic normal distributions is reasonably well for all estimators. This is more clearly illustrated in the QQ-normal plots in Figure \ref{fig4}. Using the true value of $\sigma_0^2$, as expected, the estimator $\hat{\kappa}_{n,m}$ outperforms the other methods. We compare it to our new least squares estimator in the second and fourth plots from the left. 

\begin{figure}[ht]
\centering
\includegraphics[width=11.9cm]{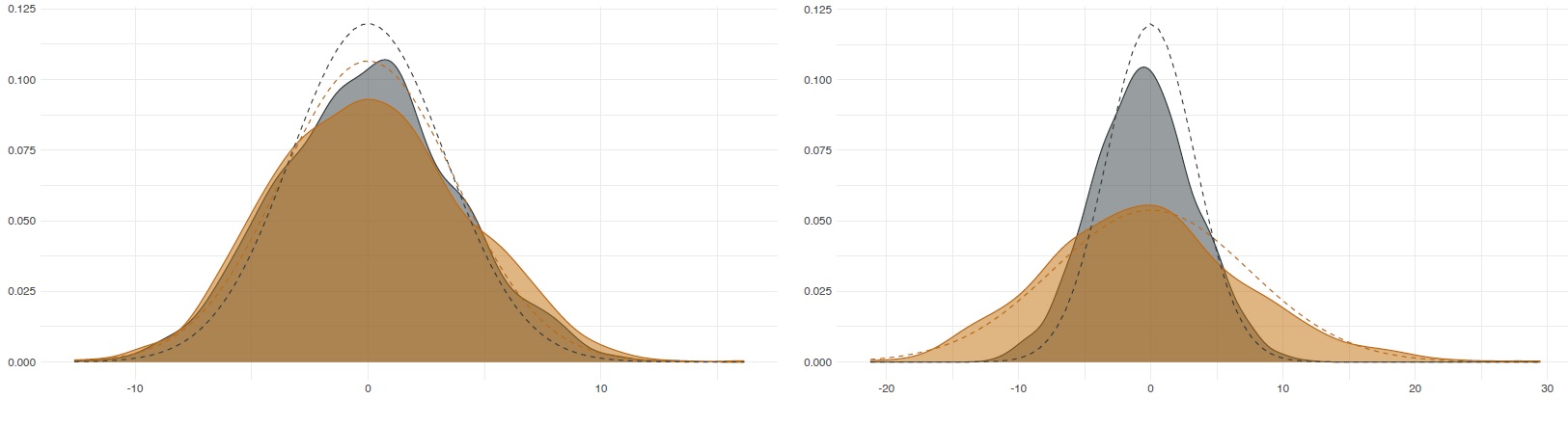}
\caption{Comparison of empirical distributions of normalized estimation errors for $\sigma_0^2$ from simulation with $n=1000$, $m=11$, $\sigma_0^2=1$, $\kappa=1$ in the left panel and $\kappa=6$ in the right panel. Grey is for $(\hat\sigma_{0}^{2})^{LS}$, and brown for the estimator by \cite{trabs}.}\label{fig3}
\end{figure}

\begin{figure}[ht]
\centering
\includegraphics[width=11.9cm]{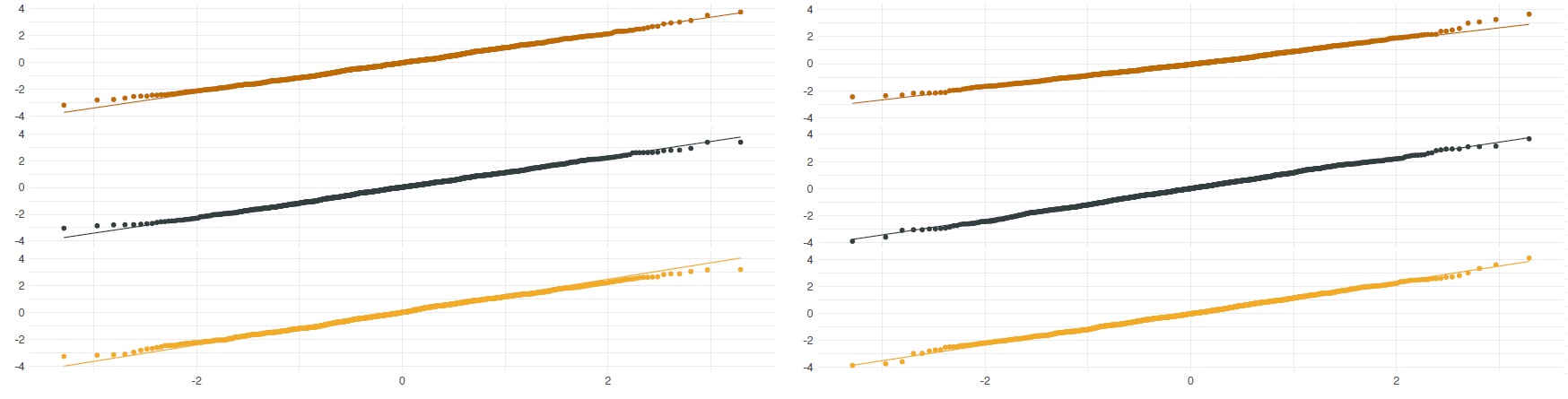}
\caption{QQ-normal plots for normalized estimation errors for $\kappa$ from simulation with $n=1000$, $m=11$, $\sigma_0^2=1$, $\kappa=1$ in the left panel and $\kappa=6$ in the right panel. Brown (top) is the estimator from \cite{trabs}, dark grey is for \eqref{lskappa} and yellow (bottom) for \eqref{kappahat}.}\label{fig4}
\end{figure}

\begin{figure}[ht]
\centering
\includegraphics[width=11.9cm]{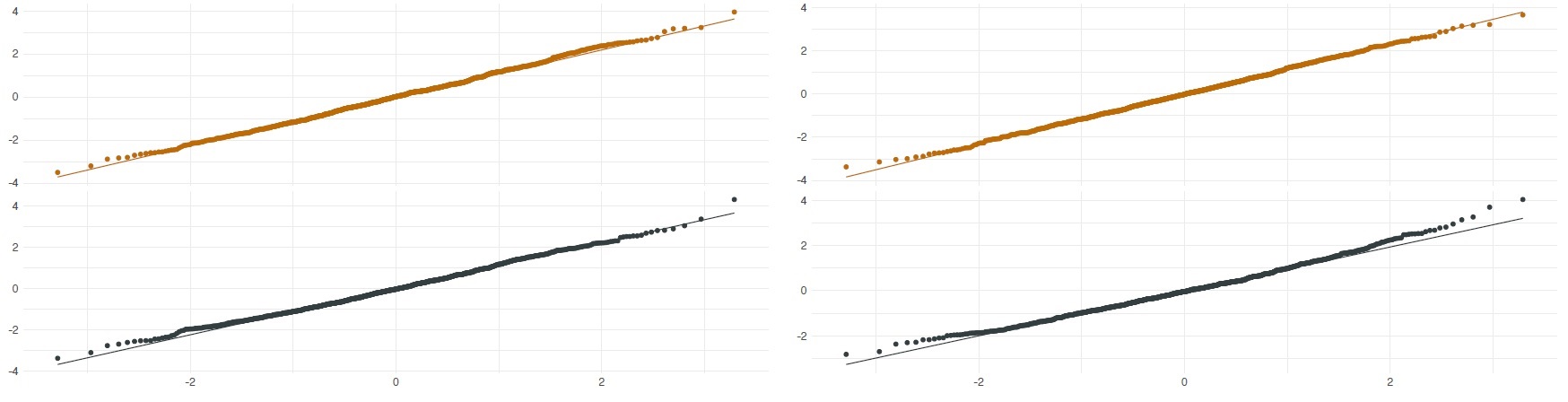}
\caption{QQ-normal plots for normalized estimation errors for $\sigma_0^2$ from simulation with $n=1000$, $m=11$, $\sigma_0^2=1$, $\kappa=1$ in the left panel and $\kappa=6$ in the right panel. Brown (top) is the estimator from \cite{trabs} and dark grey is for the estimator using \eqref{lsint}.}\label{fig5}
\end{figure}

Figure \ref{fig3} draws a similar comparison of estimated volatilities $\sigma_0^2$, for unknown $\kappa$, using the estimator \eqref{lsint} from the log-linear model and the estimator from \cite{trabs}. While for $\kappa=1$ in the left panel the performance of both methods is similar, for $\kappa=6$ in the right panel our new estimator outperforms the previous one. Figure \ref{fig5} gives the QQ-normal plots for the estimation of $\sigma_0^2$. All plots are based on 1000 Monte Carlo iterations. The QQ-normal plots compare standardized estimation errors to the standard normal. For the estimator from \cite{trabs} we use an estimated asymptotic variance based on plug-in, while for our new estimators the asymptotic variances are known constants.


\clearpage
\section{Proofs of the theorems}\label{sec5}
\subsection{Preliminaries}\label{sec5.1}
The asymptotic analysis is based on the eigendecomposition of the SPDE. The eigenfunctions $(e_k)$ and eigenvalues $(-\lambda_k)$ of the self-adjoint differential operator 
\[A_\theta =\theta_0+\theta_1\frac{\partial}{\partial y}+\theta_2\frac{\partial^{2}}{\partial y^{2}}\] 
are given by
\begin{align}\label{eq:eigenf}
  e_k(y)&=\sqrt{2}\sin\big(\pi ky\big)\exp\Big(-\frac{\theta_1}{2\theta_2}y\Big),\quad y\in [0,1],\\
  \lambda_k&=-\theta_0+\frac{\theta_1^2}{4\theta_2}+\pi^2k^2\theta_2,\qquad k\in\N\,,\label{eq:eK}
\end{align}
where all eigenvalues are negative and $(e_k)_{k\ge1}$ form an orthonormal basis of the Hilbert space $H_\theta:=\{f\colon[0,1]\to\R: \|f\|_\theta<\infty\}$ with
\begin{align*}
  \langle f,g\rangle_{\theta}&:=\int_{0}^{1}e^{y\theta_1/\theta_2}f(y)g(y)\,\d y\qquad\text{and}\qquad \|f\|_\theta^2:=\langle f,f\rangle_\theta.
\end{align*}
Let $\xi\in H_{\theta}$ be the initial condition. We impose the same mild regularity condition on $X_0=\xi$ as in Assumption 2.2 of \cite{trabs}. 
\begin{assump}\label{cond} In \eqref{eq:spde} we assume that
  \begin{enumerate}
  \item[(i)] either $\E[\langle\xi,e_k\rangle_\theta]=0$ for all $k\ge1$ and $\sup_k \lambda_k\E[\langle\xi,e_k\rangle_\theta^2]<\infty$ holds true or \\$\E[\langle A_\theta\xi,\xi\rangle_\theta]<\infty$;
  \item[(ii)] $(\langle\xi,e_k\rangle_\theta)_{k\ge1}$ are independent.
  \end{enumerate}
\end{assump}
This assumptions is more general than the one in \cite{hildebrandt} that $(X_t(y))$ is started in equilibrium and satisfied for all sufficiently regular functions $\xi$. We refer to Section 2 of \cite{trabs} for more details on the probabilistic structure.

For the solution $X_{t}(y)$ from \eqref{mild}, we have the \emph{spectral decomposition}
\begin{equation}
X_{t}(y)=\sum_{k\ge1}x_{k}(t)e_{k}(y)\,,\,\text{with}~x_{k}(t)=\langle X_{t},e_{k}\rangle_\theta\,,
\label{eq:Representation}
\end{equation}
in that the \emph{coordinate processes} $x_k$ satisfy the \emph{Ornstein-Uhlenbeck} dynamics:
\begin{equation}
\d x_{k}(t)=-\lambda_{k}x_{k}(t)\d t+\sigma_t\,\d W_{t}^{k},\quad x_{k}(0)=\langle \xi,e_{k}\rangle_{\theta}\,,\label{eq:ornsteinUhlenbeck}
\end{equation}
with independent one-dimensional Wiener processes $\{(W_{t}^{k}),k\in\N\}$.\\
We denote for some integrable random variable $Z$, its compensated version by
\begin{align}\label{comp}\overline{Z}:=Z-\E[Z]\,.\end{align}
We use upper-case characters for (compensated) random variables, while the notation for sample averages, as $\overline{y}$, is, except in Example 2 in Section 3, with lower-case characters. The notation \eqref{comp} is mainly used for $\rv(y)$ in the sequel.
\subsection{Proof of Theorem 1}\label{sec5.2}
A first-order Taylor expansion of the logarithm and Proposition 3.1 of \cite{trabs} yield that
\begin{align}
&\log\bigg(\frac{\rv(y)}{\sqrt{n}}\bigg)=\log\bigg(\E\bigg[\frac{\rv(y)}{\sqrt{n}}\bigg]+\bigg(\frac{\rv(y)}{\sqrt{n}}-\E\bigg[\frac{\rv(y)}{\sqrt{n}}\bigg]\bigg)\bigg)\notag\\
&=\log\bigg(e^{-\kappa y}\frac{\sigma_0^2}{\sqrt{\pi}}+\mathcal{O}(\Delta_n)\bigg)+\frac{\overline{\rv(y)}}{\sqrt{n}\big(e^{-\kappa y}\frac{\sigma_0^2}{\sqrt{\pi}}+\mathcal{O}(\Delta_n)\big)}+\mathcal{O}_{\P}\bigg(\bigg(\frac{\overline{\rv(y)}}{\sqrt{n}}\bigg)^2\bigg)\notag\\
&=-\kappa y+\log\bigg(\frac{\sigma_0^2}{\sqrt{\pi}}\bigg)+\mathcal{O}(\Delta_n)+\frac{\overline{\rv(y)}\sqrt{\pi}e^{\kappa y}}{\sqrt{n}\sigma_0^2}\big(1+\mathcal{O}(\Delta_n)\big)+\mathcal{O}_{\P}(\Delta_n)\notag\\
&=-\kappa y+\log\bigg(\frac{\sigma_0^2}{\sqrt{\pi}}\bigg)+\frac{\overline{\rv(y)}}{\sqrt{n}}\frac{\sqrt{\pi}e^{\kappa y}}{\sigma_0^2}+\mathcal{O}_{\P}\big(\Delta_n\big)\,,\label{taylor}
\end{align}
for some spatial point $y$. The remainders called $ R_{n,j}$ in \eqref{rvmodel}, and $\tilde R_{n,j}$ in \eqref{llm}, are contained in the last two addends. This yields for the estimator \eqref{kappahat} that
\begin{align}
\hat{\kappa}_{n,m}=\kappa-\sum_{j=1}^{m_n}\frac{\overline{\rv(y_j)}}{\sqrt{n}}\frac{y_je^{\kappa y_j} \sqrt{\pi}}{\sigma_0^2\sum_{j=1}^{m_n}y_j^2}+\mathcal{O}_{\P}(\Delta_n)\,,
\label{decomp}
\end{align}
where we conclude the order of the remainder, since under Assumption \ref{assump-obs} it holds that
\begin{align*}
\frac{\sum_{j=1}^{m_n} y_j}{\sum_{j=1}^{m_n} y_j^2}\Delta_n=\mathcal{O}(\Delta_n)\,.
\end{align*}
Since under Assumption \ref{assump-obs}, $\sqrt{nm_n}\Delta_n\to 0$, it suffices to prove a clt for the leading term from above:
\begin{align*}\sum_{i=1}^n\zeta_{n,i}:=\sqrt{m_n}\sum_{j=1}^{m_n}\overline{\rv(y_j)}\frac{\sqrt{\pi}y_je^{\kappa y_j}}{\sigma_0^2\sum_{j=1}^{m_n}y_j^2}\tod \mathcal{N}\bigg(0,\frac{3\Gamma\pi}{1-\delta+\delta^2}\bigg)\,,\end{align*} 
where $\zeta_{n,i}$ includes the $i$th squared increment of the realized volatility $\rv(y_j)$. Note that summation over time (increments) is always indexed in $i$, and summation over spatial points in $j$. Although this leading term is linear in the realized volatilities, we cannot directly adopt a clt from \cite{trabs} due to the different structure of the weights. Thus, we require an original proof of the clt for which we can reuse some ingredients from \cite{trabs}.

We begin with the asymptotic variance. We can adopt Lemma 6.4 from \cite{trabs} and Proposition 6.5 and obtain for any $\eta\in(0,1)$ that
\begin{align}
\var\bigg(\frac{1}{\sqrt{n}}\rv(y)\,e^{\kappa y}\bigg)&=\frac{\Gamma\sigma_0^4}{n}\big(1+\mathcal{O}(\Delta_n^\eta+\Delta_n^{1/2}\delta^{-1})\big)\,,\label{varrv}\\
\cov\bigg(\frac{1}{\sqrt{n}}\rv(y)\,e^{\kappa y},\frac{1}{\sqrt{n}}\rv(u)\,e^{\kappa u}\bigg)&=\mathcal{O}\bigg(\Delta_n^{3/2}\Big(\abs{y-u}^{-1}+\delta^{-1}\Big)\bigg)\,,\label{covrv}
\end{align}
for any spatial points $y$ and $u$. We obtain that
\begin{align*}
&\lim_{n\to\infty} \var\bigg(\sum_{i=1}^n\zeta_{n,i}\bigg)=\lim_{n\to\infty}\frac{m_n\pi}{\sigma_0^4\big(\sum_{j=1}^{m_n}y_j^2\big)^2}\var\bigg(\sum_{j=1}^{m_n}
\rv(y_j)\,e^{\kappa y_j}y_j\bigg)\\
&=\lim_{n\to\infty}\frac{nm_n\pi}{\sigma_0^4\big(\sum_{j=1}^{m_n}y_j^2\big)^2}\var\bigg(\frac{1}{\sqrt{n}}\sum_{j=1}^{m_n}\rv(y_j)e^{\kappa y_j}y_j\bigg)\\
&=\lim_{n\to\infty}\frac{nm_n\pi}{\sigma_0^4\big(\sum_{j=1}^{m_n}y_j^2\big)^2}\bigg(\sum_{j=1}^{m_n}y_j^2\,\var\Big(\frac{1}{\sqrt{n}}\rv(y_j)\,e^{\kappa y_j}\Big)\\
&\hspace*{3.75cm}+\sum_{j\neq  l} y_jy_l\cov\Big(\frac{1}{\sqrt{n}}\rv(y_j)\,e^{\kappa y_j},\frac{1}{\sqrt{n}}\rv(y_l)\,e^{\kappa y_l}\Big)\bigg)\\
&=\lim_{n\to\infty}\frac{nm_n\pi}{\sigma_0^4\big(\sum_{j=1}^{m_n}y_j^2\big)^2}\bigg(\frac{\Gamma\sigma_0^4\sum_{j=1}^{m_n}y_j^2}{n}\big(1+\mathcal{O}(\Delta_n^\eta)\big)\\
&\hspace*{3.75cm}+ \mathcal{O}\bigg(\Delta_n^{3/2}\Big(\sum_{j\ne  l}\frac{y_jy_l}{\abs{y_j-y_l}}+m_n^2\delta^{-1}\Big)\bigg)\bigg)\\
&=\lim_{n\to\infty} \frac{(1-2\delta)\Gamma\pi}{\frac{(1-2\delta)}{m_n}\sum_{j=1}^{m_n}y_j^2}\big(1+\mathcal{O}(\Delta_n^{\eta})\big)+\mathcal{O}\bigg(\Delta_n^{1/2}\Big(\sum_{j\ne l}\frac{y_j y_l}{m_n\abs{y_j-y_l}}+\frac{m_n}{\delta}\Big)\bigg)\\
&=\lim_{n\to\infty} \frac{(1-2\delta)\Gamma\pi}{\frac{(1-2\delta)}{m_n}\sum_{j=1}^{m_n}y_j^2}\big(1+\mathcal{O}(\Delta_n^\eta)\big)+\mathcal{O}\bigg(\Delta_n^{1/2}\Big(m_n\log(m_n)+\frac{m_n}{\delta}\Big)\bigg)\\
&=\frac{\Gamma\pi (1-2\delta)}{\int_{\delta}^{1-\delta}y^2\d y}=\frac{3\Gamma\pi (1-2\delta)}{(1-\delta)^3-\delta^3}=\frac{3\Gamma\pi }{1-\delta+\delta^2}\,.
\end{align*}
The assumption that $y_1=\delta$, $y_m=1-\delta$, is used only for the convergence of the Riemann sum in the last step. For the covariances, we used Assumption \ref{assump-obs} and an elementary estimate
\begin{align}\label{crossbound}\sum_{j\ne l}\frac{y_j y_l}{\abs{y_j-y_l}}=\mathcal{O}\bigg( \sum_{r=1}^{m_n}\sum_{l=1}^{m_n}\frac{(l+r) \, l}{m_n\,r}\bigg)=\mathcal{O}\big(m_n^2+m_n^2\log(m_n)\big)\,.\end{align}
Since $\Delta_n^{1/2}m_n\log(m_n)\to 0$ under Assumption \ref{assump-obs}, the remainders are negligible.

Next, we establish a covariance inequality for the empirical characteristic function. There exists a constant $C$, such that for all $t\in \R$:
\begin{align}\label{covineq1}
\abs{ \cov\big(\exp\big({\im tQ_a^b}\big),\exp\big({\im t Q_{b+u}^v}\big)\big)}\leq \frac{C\,t^2}{u^{3/4}}\sqrt{\var(Q_a^b)\var(Q_{b+u}^v)}\,,
\end{align}
where $Q_a^b:=\sum_{i=a}^b\zeta_{n,i}$, for natural numbers $1\leq a\leq b < b+u\leq v\leq n$.

Let $u\ge 2$, the case $u=1$ can be derived separately and is in fact easier. By the spectral decomposition \eqref{eq:Representation}
\[X_{i\Delta_n}(y_j)-X_{(i-1)\Delta_n}(y_j)=\sum_{k\ge 1}\big(x_k({i\Delta_n})-x_k({(i-1)\Delta_n})\big)e_k(y_j)\,.\]
The increments of the Ornstein-Uhlenbeck processes $(x_k(t))$ from \eqref{eq:ornsteinUhlenbeck} contain terms
\[\int_0^{(i-1)\Delta_n}e^{-\lambda_k((i-1)\Delta_n-s)}(e^{-\lambda_k \Delta_n}-1)\,\sigma \d W_s^k\,,\]
which depend on the path of $(W^k_t,0\le t\le (i-1)\Delta_n)$. Defining
\begin{align}A_2(y_j)&=\sum_{i=b+u}^{v}\bigg(\sum_{k\ge 1}\Big(x_k({i\Delta_n})-x_k({(i-1)\Delta_n})\\
&\hspace*{2cm}\notag-\int_0^{b\Delta_n} e^{-\lambda_k((i-1)\Delta_n-s)}(e^{-\lambda_k \Delta_n}-1)\,\sigma \d W_s^k\Big)e_k(y_j)\bigg)^2\,,\end{align}
we can write with the notation \eqref{comp} for squared increments
\begin{align*}
Q_{b+u}^v&=\frac{\sqrt{m_n}\sqrt{\pi}}{\sigma_0^2\sum_{j=1}^{m_n}y_j^2}\sum_{j=1}^{m_n}\sum_{i=b+u}^{v}\overline{\big(X_{i\Delta_n}-X_{(i-1)\Delta_n}\big)^2(y_j)}y_j e^{\kappa y_j}\\
&=\frac{\sqrt{m_n}\sqrt{\pi}}{\sigma_0^2\sum_{j=1}^{m_n}y_j^2}\sum_{j=1}^{m_n}\big(A_1(y_j)+A_2(y_j)\big)y_j e^{\kappa y_j}\\
&=B_1+B_2\,,\end{align*}
where $A_1$ is defined by $A_2$ and $Q_{b+u}^v$, and $B_r$, $r=1,2$, to include the sums over $A_r$. Analogous terms $A_r$ have been considered in Proposition 6.6 of \cite{trabs}. This decomposition is useful, since $B_2$ is independent of $Q_a^b$. Analogously to the proof of Proposition 6.6 in \cite{trabs}, we have for all $j$ that
\begin{align*}
\var\big({A}_1(y_j)\big)\leq \frac{\tilde C\sigma^4(v-b-u+1)\Delta_n}{(u-1)^{3/2}}\,,
\end{align*}
with some constant $\tilde C$, and from Eq.\ (59) of \cite{trabs} that
\begin{align*}
\cov\big(A_1(y_j),A_1(y_l)\big)=\mathcal{O}\bigg(\frac{\Delta_n^{3/2}(v-b-u+1)}{(u-1)^{3/2}}\frac{1}{\abs{y_j-y_l}}\bigg)\,.
\end{align*}
Thereby, we obtain that
\begin{align*}
&\var(B_1)=\frac{\pi m_n}{\sigma_0^4\big(\sum_{j=1}^{m_n}y_j^2\big)^2}\bigg(\sum_{j=1}^{m_n}e^{2\kappa y_j}y_j^2\,\var\big(A_1(y_j)\big)\\
&\hspace*{4.25cm}+\sum_{j\neq l}e^{\kappa(y_j+y_l)}y_jy_l\,\cov\big({A}_1(y_j),{A}_1(y_l)\big)\bigg)\\
&\leq \frac{\pi m_n}{\sigma_0^4\big(\sum_{j=1}^{m_n}y_j^2\big)^2}
\frac{\tilde C\sigma^4(v-b-u+1)\Delta_n}{(u-1)^{3/2}} e^{2\kappa}\sum_{j=1}^{m_n} y_j^2\\
&\hspace*{3cm}+\mathcal{O}\bigg(\frac{1}{m_n}\, \frac{\Delta_n^{3/2}(v-b-u+1)}{(u-1)^{3/2}}m_n^2\log(m_n)\bigg) \\
&\leq\frac{C'(v-b-u+1)\Delta_n}{(u-1)^{3/2}}+ \mathcal{O}\bigg(\frac{\Delta_n^{3/2}(v-b-u+1)}{(u-1)^{3/2}}m_n\log(m_n)\bigg),
\end{align*}
with a constant $C'$, where we use that $m_n\big(\sum_j y_j^2\big)^{-1}$ is bounded and \eqref{crossbound}. Since $\Delta_n^{1/2}m_n\log(m_n)\to 0$, we find a constant $C''$, such that
\begin{align}\label{upbound1}
&\var(B_1)\leq\frac{C''(v-b-u+1)\Delta_n}{(u-1)^{3/2}}\,.
\end{align}
With the variance-covariance structure of $(\zeta_{n,i})$, we obtain with some constants $C_r$, $r=1,2,3$, that
\begin{align}
\var(Q_{b+u}^{v})&\geq  C_1\frac{m_n\pi}{\sigma_0^4\big(\sum_{j=1}^{m_n}y_j^2\big)^2}\sum_{j=1}^{m_n}y_j^2\sum_{i=b+u}^v\var(\zeta_{n,i})e^{2\kappa y_j}\notag\\
&=C_2\frac{\Delta_n m_n (v-b-u+1)}{\sum_{j=1}^{m_n}y_j^2}\geq C_3 (v-b-u+1)\Delta_n\,.\label{lowbound1}
\end{align}
Since Eq.\ (54) from \cite{trabs} applied to our decomposition with $B_1$ and $B_2$, yields that
\[\abs{ \cov\big(\exp\big({\im tQ_a^b}\big),\exp\big({\im t Q_{b+u}^v}\big)\big)}\leq 2t^2\sqrt{\var(Q_a^b)\,\var(B_1)}\,,\]
\eqref{upbound1} and \eqref{lowbound1} imply \eqref{covineq1}.

A Lindeberg condition for the triangular array $(\zeta_{n,i})$ is obtained by the stronger Lyapunov condition. It is satisfied, since
\begin{align*}
&\sum_{i=1}^n\E\big[\abs{\zeta_{n,i}}^4\big]\le \frac{m_n^2\pi^2}{\sigma_0^8\big(\sum_{j=1}^{m_n}y_j^2\big)^4}\sum_{i=1}^n\E\Big[\Big(\sum_{j=1}^{m_n}\big(X_{i\Delta_n}-X_{(i-1)\Delta_n}\big)^2(y_j) y_j e^{\kappa y_j}\Big)^4\Big]\\
&\le \frac{m_n^{-2}\pi^2}{\sigma_0^8\big(m_n^{-1}\sum_{j=1}^{m_n}y_j^2\big)^4}e^{4\kappa}\sum_{i=1}^n\sum_{j,k,u,v=1}^{m_n}\E\big[\big(X_{i\Delta_n}-X_{(i-1)\Delta_n}\big)^2(y_j)\\
&\hspace*{1cm}\big(X_{i\Delta_n}-X_{(i-1)\Delta_n}\big)^2(y_k)\big(X_{i\Delta_n}-X_{(i-1)\Delta_n}\big)^2(y_u)\big(X_{i\Delta_n}-X_{(i-1)\Delta_n}\big)^2(y_v)\big]\\
&\le \frac{m_n^{-2}\pi^2}{\sigma_0^8\big(m_n^{-1}\sum_{j=1}^{m_n}y_j^2\big)^4}e^{4\kappa}\sum_{i=1}^n\sum_{j,k,u,v=1}^{m_n}\bigg(\E\big[\big(X_{i\Delta_n}-X_{(i-1)\Delta_n}\big)^8(y_j)\big]\\
&\hspace*{1cm}\times \E\big[\big(X_{i\Delta_n}-X_{(i-1)\Delta_n}\big)^8(y_k)\big]\E\big[\big(X_{i\Delta_n}-X_{(i-1)\Delta_n}\big)^8(y_u)\big]\\
&\hspace*{1cm}\times\E\big[\big(X_{i\Delta_n}-X_{(i-1)\Delta_n}\big)^8(y_v)\big]\bigg)^{1/4}\\
&=\mathcal{O}\big(m_n^2n\Delta_n^2\big)=\mathcal{O}\big(m_n^2\Delta_n\big)\,.
\end{align*}
In the last step the inner sum is estimated with a factor $m_n^4$, and we just use the regularity of $(X_t(y))_{t\ge 0}$. As $m_n^2\Delta_n\to 0$, we conclude the Lyapunov condition which together with \eqref{covineq1} and the asymptotic analysis of the variance yields the clt for the centered triangular array $(\zeta_{n,i})$ by an application of Theorem B from \cite{utev}.

\subsection{Proof of Theorem 3}\label{sec5.3}
We establish first the asymptotic variances and covariance of the estimators before proving a bivariate clt. With \eqref{taylor}, we obtain that
\begin{align*}
\hat{\kappa}_{n,m}^{LS}-\kappa=\frac{\sum_{j=1}^{m_n}\frac{\overline{\rv(y_j)}}{\sqrt{n}}\frac{\sqrt{\pi}e^{\kappa y_j}}{\sigma_0^2}\big(y_j-\overline{y}\big)}{m_n(\overline{y})^2-\sum_{u=1}^{m_n}y_u^2}+\mathcal{O}_{\P}(\Delta_n)\,,
\end{align*}
since for the remainders it holds true that
\begin{align*}
\max\bigg(\frac{\sum_{j=1}^{m_n}y_j\,\mathcal{O}_{\P}(\Delta_n)}{m_n(\overline{y})^2-\sum_{u=1}^{m_n}y_u^2},\frac{\overline{y}\,\mathcal{O}_{\P}(m_n\Delta_n)}{m_n(\overline{y})^2-\sum_{u=1}^{m_n}y_u^2}\bigg)=\mathcal{O}_{\P}(\Delta_n)\,.
\end{align*}
We can use \eqref{varrv} and \eqref{covrv} to compute the asymptotic variance:
\begin{align*}
&\lim_{n\to\infty} \var\Big(\sqrt{nm_n}(\hat{\kappa}_{n,m}^{LS}-\kappa)\Big)\\
&=\lim_{n\to\infty}\frac{nm_n \pi}{\sigma_0^4\big(m_n(\overline{y})^2-\sum_{u=1}^{m_n}y_u^2\big)^2}\var\bigg(\sum_{j=1}^{m_n}
\frac{\rv(y_j)}{\sqrt{n}}\,e^{\kappa y_j}\big(y_j-\overline{y}\big)\bigg)\\
&=\lim_{n\to\infty}\frac{nm_n\pi}{\sigma_0^4\big(m_n(\overline{y})^2-\sum_{u=1}^{m_n}y_u^2\big)^2}\bigg(\sum_{j=1}^{m_n}\big(y_j-\overline{y}\big)^2\,\var\Big(\frac{1}{\sqrt{n}}\rv(y_j)\,e^{\kappa y_j}\Big)\\
&\hspace*{1.25cm}+\sum_{j\neq  l} \big(y_j-\overline{y}\big)\big(y_l-\overline{y}\big)\cov\Big(\frac{1}{\sqrt{n}}\rv(y_j)\,e^{\kappa y_j},\frac{1}{\sqrt{n}}\rv(y_l)\,e^{\kappa y_l}\Big)\bigg)\\
&=\lim_{n\to\infty} \bigg(\frac{n m_n \pi}{\sigma_0^4\big(\sum_{u=1}^{m_n}y_u^2-m_n(\overline{y})^2\big)}\frac{\Gamma\sigma_0^4}{n}\big(1+\mathcal{O}(\Delta_n^{\eta})\big)\\
&\quad +\mathcal{O}\bigg(\frac{m_n\Delta_n^{1/2}}{\sigma_0^4\big(\sum_{u=1}^{m_n}y_u^2-m_n(\overline{y})^2\big)^2}\sum_{j\ne l}\big(y_j-\overline{y}\big) \big(y_l-\overline{y}\big)\Big(\frac{1}{\abs{y_j-y_l}}+\frac{1}{\delta}\Big)\bigg)\bigg)\\
&=\lim_{n\to\infty}\frac{\Gamma\pi}{m_n^{-1}\sum_{u=1}^{m_n}y_u^2-(\overline{y})^2}\\
&=\frac{\Gamma\pi}{(1-2\delta)^{-1}\int_{\delta}^{1-\delta}y^2\d y-\Big((1-2\delta)^{-1}\int_{\delta}^{1-\delta}y\,\d y\Big)^2}=\frac{12\Gamma\pi}{(1-2\delta)^2}\,.
\end{align*}
We used that the sum of covariances is of order
\begin{align}\label{covestimate}
\mathcal{O}\Big(m_n^{-1}\Delta_n^{1/2}\sum_{j\ne l}\frac{\big(y_j-\overline{y}\big) \big(y_l-\overline{y}\big)}{\abs{y_j-y_l}}\Big)=\mathcal{O}\Big(\Delta_n^{1/2}m_n\log(m_n)\Big)=\KLEINO(1)\,.
\end{align}
For the estimator \eqref{lsint}, we obtain that
\begin{align*}
\widehat\alpha^{LS}(\sigma_{0}^{2})=\overline{y}(\hat{\kappa}_{n,m}^{LS}-\kappa)+\alpha(\sigma_{0}^{2})+\frac{1}{m_n}\sum_{j=1}^{m_n}\frac{\overline{\rv(y_j)}}{\sqrt{n}}\frac{\sqrt{\pi}e^{\kappa y_j}}{\sigma_0^2}+\mathcal{O}_{\P}(\Delta_n),
\end{align*}
such that
\begin{align*}
\widehat\alpha^{LS}(\sigma_{0}^{2})-\alpha(\sigma_{0}^{2})=\frac{\sum_{j=1}^{m_n}\frac{\overline{\rv(y_j)}}{\sqrt{n}}\frac{\sqrt{\pi}e^{\kappa y_j}}{\sigma_0^2}\Big(y_j\overline{y}-m_n^{-1}\sum_{j=1}^{m_n}y_j^2\Big)}{m_n(\overline{y})^2-\sum_{u=1}^{m_n}y_u^2}+\mathcal{O}_{\P}(\Delta_n).
\end{align*}
With \eqref{varrv} and \eqref{covrv} and analogous steps as above, the asymptotic variance yields
\begin{align*}
&\lim_{n\to\infty} \var\Big(\sqrt{nm_n}\,\Big(\widehat\alpha^{LS}(\sigma_{0}^{2})-\alpha(\sigma_{0}^{2})\Big)\Big)\\
&=\lim_{n\to\infty}\frac{nm_n\pi}{\sigma_0^4\big(m_n(\overline{y})^2-\sum_{u=1}^{m_n}y_u^2\big)^2}\sum_{j=1}^{m_n}\Big(y_j\overline{y}-m_n^{-1}\sum_{u=1}^{m_n}y_u^2\Big)^2\,\var\Big(\frac{\rv(y_j)}{\sqrt{n}}\,e^{\kappa y_j}\Big)\\
&=\lim_{n\to\infty}\frac{m_n\Gamma\pi}{\Big(\sum_{u=1}^{m_n}y_u^2-m_n(\overline{y})^2\Big)^2}\sum_{u=1}^{m_n}y_u^2\Big(m_n^{-1}\sum_{u=1}^{m_n}y_u^2-(\overline{y})^2\Big)\\
&=\lim_{n\to\infty}\Gamma\pi\,\frac{\sum_{u=1}^{m_n}y_u^2}{\sum_{u=1}^{m_n}y_u^2-m_n(\overline{y})^2}\\
&=\Gamma\pi\,\frac{(1-2\delta)^{-1}\int_{\delta}^{1-\delta}y^2\d y}{(1-2\delta)^{-1}\int_{\delta}^{1-\delta}y^2\d y-\Big((1-2\delta)^{-1}\int_{\delta}^{1-\delta}y\,\d y\Big)^2}\,.
\end{align*}
The covariance terms for spatial points $y_j\ne y_u$ are asymptotically negligible by a similar estimate as in \eqref{covestimate}. The asymptotic covariance between both estimators yields
\begin{align*}
&\lim_{n\to\infty} nm_n\cov\Big(\widehat\alpha^{LS}(\sigma_{0}^{2})\,,\,\hat{\kappa}_{n,m}^{LS}\Big)\\
&=\lim_{n\to\infty}\frac{nm_n\pi}{\sigma_0^4\big(m_n(\overline{y})^2-\sum_{u=1}^{m_n}y_u^2\big)^2}\sum_{j=1}^{m_n}\Big(y_j\overline{y}-\frac{\sum_{u=1}^{m_n}y_u^2}{m_n}\Big)\Big(y_j-\overline{y}\Big)\,\var\Big(\frac{\rv(y_j)}{\sqrt{n}}\,e^{\kappa y_j}\Big)\\
&=\lim_{n\to\infty}m_n\Gamma\pi\,\frac{\overline{y}\Big(\sum_{u=1}^{m_n}y_u^2-m_n(\overline{y})^2\Big)}{\Big(\sum_{u=1}^{m_n}y_u^2-m_n(\overline{y})^2\Big)^2}\\
&=\lim_{n\to\infty}\Gamma\pi\,\frac{\overline{y}}{m_n^{-1}\sum_{u=1}^{m_n}y_u^2-(\overline{y})^2}\\
&=\Gamma\pi\,\frac{(1-2\delta)^{-1}\int_{\delta}^{1-\delta}y\,\d y}{(1-2\delta)^{-1}\int_{\delta}^{1-\delta}y^2\d y-\Big((1-2\delta)^{-1}\int_{\delta}^{1-\delta}y\,\d y\Big)^2}\,.
\end{align*}
The covariance terms for spatial points $y_j\ne y_u$ are asymptotically negligible by a similar estimate as in \eqref{covestimate}. 
Computing the elementary integrals and simple transformations yield the asymptotic variance-covariance matrix $\Sigma$ in Theorem \ref{jointclt}.

To establish the bivariate clt, it suffices to prove the clt for the $\R^2$-valued triangular array
\begin{align*}\Xi_{n,i}\hspace*{-.05cm}=\hspace*{-.05cm}\frac{\sqrt{m_n\pi}}{\sigma_0^2\big(m_n(\overline{y})^2-\sum_{u=1}^{m_n}y_u^2\big)}\sum_{j=1}^{m_n}\overline{\big(X_{i\Delta_n}\hspace*{-.075cm}-\hspace*{-.05cm}X_{(i-1)\Delta_n}\big)^2(y_j)}e^{\kappa y_j}\hspace*{-.075cm}\left(\begin{array}{c}y_j-\overline{y}\\[.05cm] y_j\,\overline{y}-\frac{\sum_{u=1}^{m_n}y_u^2}{{m_n}}\end{array}\right).\end{align*}
Here, we use the notation \eqref{comp} for the squared time increments. The first entry of this vector is the leading term of \(\sqrt{nm_n}(\hat{\kappa}_{n,m}^{LS}-\kappa)\), and the second entry of \(\sqrt{nm_n}\,\big(\widehat\alpha^{LS}(\sigma_{0}^{2})-\alpha(\sigma_{0}^{2})\big)\). We apply the Cram\'{e}r-Wold device and Theorem B from \cite{utev}. Taking the scalar product with some arbitrary $\gamma\in\R^2$, we obtain by linearity that 
\begin{align*}\langle \gamma,\Xi_{n,i}\rangle &=S_{mn}\sum_{j=1}^{m_n}\overline{\big(X_{i\Delta_n}-X_{(i-1)\Delta_n}\big)^2(y_j)}e^{\kappa y_j} G_j^{\gamma}\,,\\
\mbox{with}~~~~S_{mn}&=\frac{\sqrt{m_n\pi}}{\sigma_0^2\big(m_n(\overline{y})^2-\sum_{u=1}^{m_n}y_u^2\big)}\,,\\[.075cm]
\mbox{and}~~~~G_j^{\gamma}&=\bigg\langle \gamma, \left(\begin{array}{c}y_j-\overline{y}\\ y_j\,\overline{y}-\frac{1}{m_n}\sum_{u=1}^{m_n}y_u^2\end{array}\right)\bigg\rangle\,.
\end{align*}
Note that for any $\gamma\in\R^2$, $S_{mn}\,G_j^{\gamma}$ is uniformly in $j$ bounded by a constant, such that the structure for proving a covariance inequality for the empirical characteristic function and a Lyapunov condition is analogous to the one-dimensional case. With $\xi_{n,i}^{\gamma}:=\langle \gamma,\Xi_{n,i}\rangle$, 
\[\sum_{i=1}^n\xi_{n,i}^{\gamma}=\Big\langle \gamma,\sum_{i=1}^n\Xi_{n,i}\Big\rangle=\sum_{i=1}^{n}\langle \gamma,\Xi_{n,i}\rangle\,,\] 
we obtain for $\tilde Q_a^b:=\sum_{i=a}^b\xi_{n,i}^{\gamma}$, that
\begin{align*}
\tilde Q_{b+u}^v&=S_{mn}\sum_{j=1}^{m_n}\big(A_1(y_j)+A_2(y_j)\big)y_j e^{\kappa y_j}\,G_j^{\gamma},\end{align*}
with the same terms $A_1(y_j)$ and $A_2(y_j)$ as in the proof of \eqref{covineq1}. Therefore, using the same bounds as in the proof of \eqref{covineq1}, we obtain that
\begin{align}\label{covineq2}
\abs{ \cov\big(\exp\big({\im t\tilde Q_a^b}\big),\exp\big({\im t \tilde Q_{b+u}^v}\big)\big)}\leq \frac{C\,t^2}{u^{3/4}}\sqrt{\var\big(\tilde Q_a^b\big)\var\big(\tilde Q_{b+u}^v\big)}\,,
\end{align}
for all $t\in \R$, for natural numbers $1\leq a\leq b < b+u\leq v\leq n$, and for some constant $C$.

The Lyapunov condition for the triangular array $(\xi_{n,i}^{\gamma})$ holds, since
\begin{align*}
&\sum_{i=1}^n\E\big[\abs{\xi_{n,i}^{\gamma}}^4\big]=\sum_{i=1}^n S_{mn}^4 \E\Big[\Big(\sum_{j=1}^{m_n}\overline{\big(X_{i\Delta_n}-X_{(i-1)\Delta_n}\big)^2(y_j)}  e^{\kappa y_j}G_j^{\gamma}\Big)^4\Big]\\
&\le C\,e^{4\kappa}\sum_{i=1}^n\sum_{j,k,u,v=1}^{m_n}\E\big[\big(X_{i\Delta_n}-X_{(i-1)\Delta_n}\big)^2(y_j)\big(X_{i\Delta_n}-X_{(i-1)\Delta_n}\big)^2(y_k)\\
&\hspace*{1cm}\big(X_{i\Delta_n}-X_{(i-1)\Delta_n}\big)^2(y_u)\big(X_{i\Delta_n}-X_{(i-1)\Delta_n}\big)^2(y_v)\big]\\
&\le C\,e^{4\kappa}\sum_{i=1}^n\sum_{j,k,u,v=1}^{m_n}\bigg(\E\big[\big(X_{i\Delta_n}-X_{(i-1)\Delta_n}\big)^8(y_j)\big]\\
&\hspace*{1cm}\times \E\big[\big(X_{i\Delta_n}-X_{(i-1)\Delta_n}\big)^8(y_k)\big]\E\big[\big(X_{i\Delta_n}-X_{(i-1)\Delta_n}\big)^8(y_u)\big]\\
&\hspace*{1cm}\times\E\big[\big(X_{i\Delta_n}-X_{(i-1)\Delta_n}\big)^8(y_v)\big]\bigg)^{1/4}\\
&=\mathcal{O}\big(m_n^2n\Delta_n^2\big)=\mathcal{O}\big(m_n^2\Delta_n\big)\,.
\end{align*}
for some constant $C$. As $m_n^2\Delta_n\to 0$, we conclude the Lyapunov condition which together with \eqref{covineq2} and the asymptotic variance-covariance structure yields the clt for the triangular array $(\xi_{n,i}^{\gamma})$, for any $\gamma\in\R^2$, by an application of Theorem B from \cite{utev}. We conclude with the Cram\'{e}r-Wold device.
\subsection{Proof of Theorem 2}\label{sec5.4}
Theorem 2 is established as a simple corollary of Theorem 3 showing that the two estimators \eqref{kappahatf} and \eqref{lskappa} coincide. This is based on the formula that for vectors $y,z\in\mathds{R}^m$, we have that
\begin{align}\label{esteq}\sum_{j\ne l}(z_j-z_l)(y_j-y_l)=2\,m\,\sum_{j=1}^m\big(z_j-\overline{z}\big)\big(y_j-\overline{y}\big)=2\,m\,\sum_{j=1}^m z_j\big(y_j-\overline{y}\big)\,,\end{align}
using our standard notation for means $\overline{y}$ and $\overline{z}$ applied to the vectors. \eqref{esteq} is true, since
\begin{align*}
\sum_{j\ne l}(z_j-z_l)(y_j-y_l)&=\sum_{j, l=1}^m(z_j-z_l)(y_j-y_l)\\
&=m\sum_{j=1}^m z_j\,y_j-2\,\sum_{j, l=1}^my_j\, z_l+m\sum_{l=1}^m z_l \,y_l\\
&=2m \sum_{j=1}^m z_j\,y_j-2m^2\,\overline{y}\,\overline{z}\,,
\end{align*}
and by the transformation
\begin{align*}
\sum_{j=1}^m\big(z_j-\overline{z}\big)\big(y_j-\overline{y}\big)=\sum_{j=1}^m z_j\,y_j - m\,\overline{y}\,\overline{z}\,.
\end{align*}
Applying \eqref{esteq} twice, to the numerator and to the denominator of \eqref{kappahatf} yields the estimator \eqref{lskappa}. We hence conclude the clt in Theorem 2 as the marginal clt from the bivariate clt given in Theorem 3.
\section*{Acknowledgments}
The authors are grateful to two anonymous reviewers and Florian Hildebrandt for their helpful comments.
%
%
%

\section*{Declarations}
Both authors declare that they have no conflicts of interest to disclose.
\addcontentsline{toc}{section}{References}
\bibliography{library}


\end{document}